\documentclass[12pt]{article}
\usepackage{amsfonts, amsmath}
\newcommand{\ggg}{\mathbb G}

\newcommand{\SSS}{{\cal S}}

\newtheorem{theorem}{Theorem}
\newtheorem{rem}{Remark}
\newtheorem{defn}{Definition}

\author{Mark A. Pankov}
\title{Sets of singular restrictions of
symplectic forms}
\date{Institute of Mathematics,
Kiev, Ukraine\\
e-mail:\;\texttt{pankov@imath.kiev.ua}}

\begin{document}

\maketitle

\section{Introduction}

Let $V$ be an $n$-dimensional vector space over
some field $F$.
We start with the following well-known
definitions.

\begin{defn}{\rm
For an automorphism $\sigma \in Aut(F)$ and
other vector space
$V'$ over the field $F$ a mapping $f:V\to V'$ is called
$\sigma$-{\it linear} if
$$f(x+y)=f(x)+f(y)\;\;\;\;\forall\;x,y\in V$$
and
$$f(ax)=\sigma(a)f(x)\;\;\;\;\forall\;x\in V, \;a\in F\;.$$
In the case when
$$Ker f={0}\;\mbox{ and }\;f(V)=V'$$
$f$ is a bijection of $V$ onto $V'$ and
we say that $f$ is a $\sigma$-{\it linear isomorphism}
of $V$ onto $V'$.
Moreover, if $V=V'$ then the isomorphism $f$ will be
called a $\sigma$-{\it linear transformation} of $V$.
We shall often exploited the term "linear" in plase of
$\sigma$-linear.}
\end{defn}

\begin{defn}{\rm
Let $\sigma_{1},\sigma_{2}\in Aut(F)$.
A mapping
$$\Omega:V\times V\to F$$
is a $(\sigma_{1},\sigma_{2})$-{\it bilinear
form} on $V$ if
$$\Omega(\cdot , x)\;\mbox{ and }\;\Omega(x, \cdot)$$
are $\sigma_{1}$-linear and  $\sigma_{2}$-linear
functionals on $V$.
In what follows we shall often say that a
form is "bilinear";
the full name will be used only for some special
cases.
}\end{defn}

\begin{defn}{\rm
A bilinear form $\Omega$ defined on
$V$ is called {\it symplectic} if the equality $\Omega(x,x)=0$
holds for each $x\in V$.
}\end{defn}

\begin{rem}{\rm
Recall that a bilinear form $\Omega$ on $V$
is called {\it skew-symmetric} if
$$\Omega(x,y)=-\Omega(y,x)\;\;\;\;\forall\;x,y\in V\;.$$
Each symplectic form is skew-symmetric.
The inverse statement holds only if the
characterystic of the field $F$ is not equal to $2$.
In the case
when $char F=2$ there exists a skew-symmetric forn which is not
symplectic.}
\end{rem}

\begin{rem}{\rm
An immediate verification shows that for each
skew-symmetric $(\sigma_{1},\sigma_{2})$-bilinear form $\Omega$
we have $\sigma_{1}=\sigma_{2}$.
}\end{rem}

Now assume that our symplectic form $\Omega$ is
{\it non-singular} (then the number $n$ is even).
Denote by $\ggg^{n}_{k}$ the Grassmannian
manifold of $k$-dimensional planes contained in $V$
and consider the set
$S^{n}_{k}(\Omega)$ of all
planes $s\in \ggg^{n}_{k}$ such that
the restriction of the form $\Omega$ onto $s$
is singular.

\begin{rem}{\rm
The forms $a \Omega$ and $\sigma \Omega$
are symplectic for each $a \in F$ and $\sigma\in Aut(F)$.
Moreover,
$$S^{n}_{k}(\Omega)=S^{n}_{k}(a\Omega)=S^{n}_{k}(\sigma\Omega)\;.$$
}\end{rem}

Any two non-singular symplectic form
$\Omega_{1}$ and $\Omega_{2}$ on $V$
are {\it equivalent}; i.e. there exists a
$\sigma$-linear transformation $f$ of $V$ such that
$$\Omega_{2}(x,y)=\Omega_{1}(f(x),f(y))
\;\;\;\forall\;x,y\in V\;.$$
This implies that the set
$S^{n}_{k}(\Omega_{2})$ can be transfered
to $S^{n}_{k}(\Omega_{1})$ by the transformation
of $\ggg^{n}_{k}$ induced by $f$.

\begin{rem}{\rm
Each $\sigma$-linear transformation of $V$
($\sigma \in Aut(V)$)
induces some transformation of
the Grassmannian manifold $\ggg^{n}_{k}$.
This fact will be often exploited in what follows.
}\end{rem}

Any symplectic form defined on an odd-dimensional
vector space is singular.
Therefore, the sets $S^{n}_{k}(\Omega)$
and $\ggg^{n}_{k}$ are coincident
if the number $k$ is odd.
For the case when $k$ is even the structure
of $S^{n}_{k}(\Omega)$ is more
complicated.
However, for the cases $k=2,n-2$
there exists a simple geometrical characteristic
of this set. We shall study it above.

Let $s\in \ggg^{n}_{m}$ and $m\ne k$.
Denote by $\ggg^{n}_{k}(s)$
the set of all planes $l\in\ggg^{n}_{k}$
satisfying the condition
$$l\subset s \mbox{ if } m>k$$
or
$$s\subset l \mbox{ if } m<k\;.$$
It must be pointed out that $\ggg^{n}_{k}(s)$
can be
considered as some Grassmannian manifold.

\begin{defn}{\rm
Let $X\subset \ggg^{n}_{k}$ and $k=2$ or $n-2$.
In the case $k=n-2$ we say that $X$ satisfies
condition $\SSS$ if for each plane $s\in \ggg^{n}_{n-1}$
there exists a line $F(s)$ contained in $s$
and such that
\begin{equation}
\ggg^{n}_{k}(s)\cap X=\ggg^{n}_{k}(s)\cap\ggg^{n}_{k}(F(s))\;;
\end{equation}
in other words, the set $\ggg^{n}_{k}(s)\cap X$
consists of all planes containing some line.
In the case $k=2$ the set $X$ satisfies
condition $\SSS$ if for each line $s\in \ggg^{n}_{1}$
there exists a plane $F(s)\in\ggg^{n}_{n-1}$
containing $s$ and such that
equality (1) holds true.}
\end{defn}

It is not difficult to prove (Section 3)
that   the set
$S^{n}_{k}(\Omega)$ ($k=2,n-2$) satisfies  condition
$\SSS$.
This paper is devoted to prove the inverse
statement.

\begin{theorem}
For a set $X \subset \ggg^{n}_{k}$,
$k=2$ or $n-2$, the following two conditions
are equivalent:
\begin{enumerate}
\item[\rm{(i)}] $X$ satisfies condition $\SSS$,
\item[\rm{(ii)}] there exists a non-singular
symplectic form $\Omega$ on $V$ such that
$X=S^{n}_{k}(\Omega)$.
\end{enumerate}
\end{theorem}

For the case when $k\ne 2,n-2$ (clearly, the number $k$ is even)
the similar characteristic of the set $S^{n}_{k}(\Omega)$
is  not found yet.

\section{Bijections of $\ggg^{n}_{k}$ onto
$\ggg^{n}_{n-k}$ defined by bilinear forms}

Let $\Omega$ be a non-singular bilinear form on
$V$ and $s\in \ggg^{n}_{k}$.
Denote by $s^{\perp}_{\Omega}$ the $\Omega$-orthogonal
complement to $s$; i.e.
$$s^{\perp}_{\Omega}=
\{\;y\in V\;|\;\Omega(x,y)=0\;\;\;\forall\;x\in s\;\}\;.$$
Then $s^{\perp}_{\Omega}\in \ggg^{n}_{n-k}$
and the form $\Omega$
defines the
mapping of
$$F^{n}_{k\,n-k}(\Omega):\ggg^{n}_{k}\to\ggg^{n}_{n-k}$$
tranferring each plane
to the $\Omega$-orthogonal complement.
It is a bijection onto.

We shall use the following properties
of bijections defined by non-singular bilinear forms.
\begin{enumerate}
\item[A] {\it For each non-singular form $\Omega$
on $V$ there exists a non-singular form
$\Omega'$ such that }
$$(F^{n}_{k\,n-k}(\Omega))^{-1}=F^{n}_{n-k\,k}(\Omega')\;.$$
Define
$$\Omega'(x,y)=\Omega(y,x)\;\;\;\;\forall\;x,y\in V\;.$$
It is trivial that the form $\Omega'$ is non-singular
and satisfies  the required condition.

\item[B] {\it The equality
$$(F^{n}_{k\,n-k}(\Omega))^{-1}=F^{n}_{n-k\,k}(\Omega)$$
holds for some non-singular bilinear forn
$\Omega$
if and only if this form is reflexive};
i.e. the following two conditions
$$\Omega(x,y)=0\;\mbox{ and }\;\Omega(y,x)=0$$
are equivalent for any $x,y\in V$.

\item[C] {\it For  any two non-singular
bilinear forms $\Omega_{1}$ and $\Omega_{2}$ on $V$
there exist transformations $f$ and $g$ of
$\ggg^{n}_{k}$ and $\ggg^{n}_{n-k}$
{\rm(}respectively{\rm)} induced by linear transformations of
$V$ and such that}
$$F^{n}_{k\,n-k}(\Omega_{2})=
F^{n}_{k\,n-k}(\Omega_{1})f=gF^{n}_{k\,n-k}(\Omega_{1})$$
(see [4]).

\item[D] Properties A and C imply that
{\it for  any two non-singular
bilinear forms $\Omega_{1}$ and $\Omega_{2}$ on $V$
the composition
$$F^{n}_{n-k\,k}(\Omega_{2})F^{n}_{k\,n-k}(\Omega_{1})$$
is a transformation of $\ggg^{n}_{k}$ induced by
some $\sigma$-linear
transformation of $V$}.
\end{enumerate}

\section{Proof of implication $\rm{(ii)}\Rightarrow \rm{(i)}$}

Let $\Omega$ be a non-singular symplectic form
on $V$.
It is trivial that a planes  $l\in \ggg^{n}_{k}$
belongs to the set $S^{n}_{k}(\Omega)$ if and only if
\begin{equation}
\dim l\cap l^{\perp}_{\Omega} \ge 1\;.
\end{equation}
This implies the following inclusion
$$
\ggg^{n}_{n-2}(s)\cap\ggg^{n}_{n-2}(s^{\perp}_{\Omega})
\subset
\ggg^{n}_{n-2}(s)\cap S^{n}_{n-2}(\Omega)
\;\;\;\;\forall\;s\in \ggg^{n}_{n-1}\;.$$
Assume that the inverse inclusion fails and there exists
$$l\in\ggg^{n}_{n-2}(s)\cap S^{n}_{n-2}(\Omega)$$
which does not contain the line $s^{\perp}_{\Omega}$.
The condition $l\in S^{n}_{n-2}(\Omega)$
guarantees the existence of a non-zero vector $x\in l$
such that $\Omega(x,y)=0$ for each vector
$y\in l$. The plane $s$ is generated by
the plane $l$ and the line $s^{\perp}_{\Omega}$.
Therefore, the similar equality
holds for each vector $y\in s$; i.e.
$x\in s^{\perp}_{\Omega}$. These arguments
dispruve our hypothesis and the requred inclusion is
proved. In other words,
we ontain the equality
\begin{equation}
\ggg^{n}_{n-2}(s)\cap S^{n}_{n-2}(\Omega)=
\ggg^{n}_{n-2}(s)\cap\ggg^{n}_{n-2}(s^{\perp}_{\Omega})
\;\;\;\;\forall\;s\in \ggg^{n}_{n-1}
\end{equation}
showing that $S^{n}_{n-2}(\Omega)$ satisfies  condition
$\SSS$.

It was noted above that a plane $l\in \ggg^{n}_{k}$
belongs to the set $S^{n}_{k}(\Omega)$ if and only if
equation (2) holds. This implies that
the bijection $F^{n}_{k\,n-k}(\Omega)$
transfers $S^{n}_{k}(\Omega)$ to $S^{n}_{n-k}(\Omega)$.
Moreover, $F^{n}_{n-2\,2}(\Omega)$ transfers equation (3)
to the following equality
$$
\ggg^{n}_{2}(s)\cap S^{n}_{2}(\Omega)=
\ggg^{n}_{2}(s)\cap\ggg^{n}_{2}(s^{\perp}_{\Omega})
\;\;\;\;\forall\;s\in \ggg^{n}_{1}\;;
$$
i.e. $S^{n}_{2}(\Omega)$ satisfies  condition
$\SSS$.

\section{Reduction of the case $k=2$ to the case $k=n-2$}

Assume that Theorem 1 is proved for the case $k=n-2$
and consiner a set $X \subset \ggg^{n}_{2}$
satisfying condition $\SSS$.
Let $g$ be a
bijection of $\ggg^{n}_{2}$ onto
$\ggg^{n}_{n-2}$ defined by some
bilinear form. An immediate verification shows that
$g(X)$ is a subset of $\ggg^{n}_{n-2}$
satisfying  condition $\SSS$ and
there exists a non-singular symplectic form
$\Omega$ such that
$$g(X)=S^{n}_{n-2}(\Omega)\;.$$
Then
$$F^{n}_{n-2\,2}(\Omega)g(X)=S^{n}_{2}(\Omega)\;.$$
It was noted above that
$F^{n}_{n-2\,2}(\Omega)g$
is a transfomation
of $\ggg^{n}_{2}$ induced by a $\sigma$-linear
trandformation $f$ of $V$ (property D).
This implies that
$$
X=S^{n}_{2}(\Omega')\;,
$$
where $\Omega'$ is the non-singular symplectic form
defined by the condition
$$\Omega'(x,y)=\Omega(f(x),f(y))\;\;\;\;
\forall\;x,y\in V\;.$$

\section{Proof of implication $\rm{(i)}\Rightarrow \rm{(ii)}$}

Each set $X \subset \ggg^{n}_{n-2}$
satisfying condition $\SSS$
is defined by some mapping
$$F:\ggg^{n}_{n-1}\to\ggg^{n}_{1}\;.$$

\begin{rem}{\rm
It is trivial that
for the case when $X=S^{n}_{n-2}(\Omega)$
the mapping $F$ coincides with
$F^{n}_{n-1\,1}(\Omega)$.
}\end{rem}

We state that the mapping $F$
satisfies the
following two conditions:
\begin{enumerate}
\item[{\rm (1)}] each plane
$s\in \ggg^{n}_{n-1}$ contains the line $F(s)$;
\item[{\rm (2)}] for any two planes
$s_{1},s_{2}\in \ggg^{n}_{n-1}$
the line $F(s_{1})$ is contained in $s_{2}$
if and only if $F(s_{2})$ is contained in $s_{1}$.
\end{enumerate}

First condition is trivial.

To prove second condition consider the plane
$$l=s_{1}\cap s_{2}\in \ggg^{n}_{n-2}\;.$$
If $l\notin X$ then
the lines $F(s_{1})$ and $F(s_{2})$
is not contained in $l$.
Then condition (1) guarantees that
$F(s_{1})$ is not contained in $s_{2}$
and
$F(s_{2})$ is not contained in $s_{1}$.
In the case when $l\in X$ the plane
$l$ contains $F(s_{1})$ and $F(s_{2})$.
The fulfilment of condition (2) is proved.

Let $\Omega$ be a non-singular bilinear form on $V$.
Recall that a vector $y\in V$ is called $\Omega$-{\it orthogonal}
to a vector $x\in V$ if $\Omega(x,y)=0$.  If our form is
reflexive then $y$ is $\Omega$-orthogonal to $x$ if and only if
$x$ is $\Omega$-orthogonal to $y$.  In this case we can say that
the vectors $x$ and $y$ are $\Omega$-{\it orthogonal}.

Now fix some reflexive non-singular bilinear form $\Omega$ on $V$
and consider the mapping
$$f=FF^{n}_{1\,n-1}(\Omega):\ggg^{n}_{1}\to \ggg^{n}_{1}\;.$$
It is not difficult to see that it
satisfies the next conditions:
\begin{enumerate} \item[(1')]
for each $l\in \ggg^{n}_{1}$ the lines $l$ and
$f(l)$ are $\Omega$-orthogonal;
\item[(2')]
for any two lines $t_{1},t_{2}\in \ggg^{n}_{1}$ the line $t_{1}$
and $f(t_{2})$ are $\Omega$-orthogonal if and only if $t_{2}$ and
$f(t_{1})$ are $\Omega$-orthogonal.
\end{enumerate}
We want to show that {\it the  mapping $f$ is
induced by some $\sigma$-linear transformation
of $V$}.
Then
$$F=f(F^{n}_{1\,n-1}(\Omega))^{-1}=fF^{n}_{n-1\,1}(\Omega)$$
(the form $\Omega$ is reflexive) and
property C implies
the existence of a non-singular bilinear form
$\Omega'$ on $V$ such that
$$F=F^{n}_{n-1\,1}(\Omega')\;.$$
Condition (1) guarantees that
any plane $s\in \ggg^{n}_{n-1}$
contains the $\Omega'$-orthogonal complement
$s^{\perp}_{\Omega'}$. It is not difficult to see that
each form satisfying this condition is symplectic
and the required statement will be proved.

First of all show that {\it the mapping $f$
is surjective.}
Let $l \in \ggg^{n}_{1}$ and
$l_{1},...,l_{n-1}$
be linearly independent lines generating the
plane $l^{\perp}_{\Omega}$.
Then the lines $f(l_{1}),...,f(l_{n-1})$
generate a plane $s$ the dimension of which is not greater
than $n-1$.
Consider a line $l'$ contained in $s^{\perp}_{\Omega}$.
Condition (2') shows that the line $f(l')$ is
$\Omega$-orthogonal to the lines $l_{1},...,l_{n-1}$;
therefore, $f(l')=l$.

We state that
{\it for each plane
$s\in \ggg^{n}_{n-1}$ there exist
planes $s_{1}, s_{2}\in \ggg^{n}_{n-1}$
such that
$$
f(\ggg^{n}_{1}(s))=\ggg^{n}_{1}(s_{1})
$$
and}
$$
f^{-1}(\ggg^{n}_{1}(s))=\ggg^{n}_{1}(s_{2})\;.
$$
Let $l\in \ggg^{n}_{1}$ be the $\Omega$-orthogonal complement to
the plane $s$.
The mapping $f$ is a surjection and
there exists a line $l_{1} \in \ggg^{n}_{1}$ such that
$f(l_{1})=l$.
Denote by $s_{1}$ the
$\Omega$-orthogonal complement to
$l_{1}$. After that denote by $s_{2}$ the orthogonal complement to
the line $f(l)$.
Then condition (2')
implies the fulfilment of the requred equalities.

Now assume that
for some line $l\in \ggg^{n}_{1}$ there exist
two lines $l_{1},l_{2}\in \ggg^{n}_{1}$
such that
$$f(l_{i})=l\;\;\;\;i=1,2\;.$$
The previous arguments show that
$$f(\ggg^{n}_{1}(l^{\perp}_{\Omega}))=
\ggg^{n}_{1}((l_{1})^{\perp}_{\Omega})$$
and
$$f(\ggg^{n}_{1}(l^{\perp}_{\Omega}))=
\ggg^{n}_{1}((l_{2})^{\perp}_{\Omega})\;.$$
Then the planes
$(l_{1})^{\perp}_{\Omega}$
and
$(l_{2})^{\perp}_{\Omega}$
are coincident; i.e.
$l_{1}=l_{2}$.
In other words, we have proved
that {\it the mapping $f$ is a injective}.

Therefore, $f$ is a bijection onto.
Moreover, $f$ and $f^{-1}$ transfer each
collection of linearly independent lines to a
collection of linearly independent lines.
The Fundamental Theorem of Projective Geometry
[1--3]
states that in this case
$f$ is induced by some $\sigma$-linear transformation
of $V$.

\begin{enumerate}
\item[[1]] E. Artin, {\it Geometric Algebra}, Interscience,
New York, 1957.

\item[[2]] J. Dieudonne,
{\it La Geometrie des groupes classiques}, Springer -- Verlag,
Berlin -- New York, 1971.

\item[[3]] O. T. O'Meara, {\it Lectures on linear groups},
Providence, Rhode Island, 1974.

\item[[4]] M. A. Pankov, {\it Irregular subsets of the
Grassmannian manifolds}, 1999 (book, math.AT/9910081).
\end{enumerate}

\end{document}